\documentclass[11pt]{article}
\usepackage{amscd}
\usepackage{amsfonts}
\usepackage{amsmath}
\usepackage{amssymb}
\usepackage{amsthm}
\usepackage{bbm}
\usepackage{CJK}
\usepackage{fancyhdr}
\usepackage{graphicx}
\usepackage{hyperref}
\usepackage{indentfirst}
\usepackage{latexsym}
\usepackage{mathrsfs}
\usepackage{xypic}

\usepackage[top=1in,bottom=1in,left=1.25in,right=1.25in]{geometry}
\def\<{\langle}
\def\>{\rangle}
\def\a{\alpha}
\def\b{\beta}

\def\c{\cdot}

\def\D{\Delta}

\def\g{\gamma}

\def\o{\otimes}

\def\r{\rho}

\def\v{\varepsilon}

\date{}
\begin{document}
\renewcommand{\baselinestretch}{1.2}
\renewcommand{\arraystretch}{1.0}
\title{\bf{ Relative Hom-Hopf modules and total integrals }}
\author{{\bf Shuangjian Guo$^{1}$, Xiaohui Zhang$^{2}$
        and Shengxiang Wang$^{3}$\footnote
        {Correspondence: S X Wang: wangsx-math@163.com.} }\\
1.~ School of Mathematics and Statistics, Guizhou University of Finance and\\ Economics, Guizhou  550025, China\\
2.~Department of Mathematics, Southeast  University,
Nanjing 210096, China\\
 3.~ School of Mathematics and Finance, Chuzhou University,
 Chuzhou 239000,  China}
 \maketitle

\noindent\textbf {\bf Abstract} Let $(H, \a)$ be a monoidal Hom-Hopf algebra and $(A, \b)$ a right $(H, \a)$-Hom-comodule algebra.
We first investigate the criterion for the existence
of a total integral of $(A, \b)$ in the setting of monoidal Hom-Hopf algebras.
Also we  prove that there exists a total integral $\phi: (H, \a)\rightarrow (A, \b)$ if and only if any representation
of the pair $(H,A)$ is injective in a functorial way, as a corepresentation of $(H, \a)$,
which generalizes Doi's result. Finally, we define a total quantum integral $\g: H\rightarrow Hom(H, A)$ and prove the following affineness
criterion: if there exists  a total quantum integral $\g$ and the
canonical map $
\psi: A\o_{B}A\rightarrow A\o H,\ \
a\o_{B}b\mapsto \b^{-1}(a)b_{[0]}\o \a(b_{[1]})
$
is surjective, then the induction functor
$A\o_B-: \widetilde{\mathscr{H}}(\mathscr{M}_k)_{B}\rightarrow \widetilde{\mathscr{H}}(\mathscr{M}_k)^{H}_{A}$
is an equivalence of categories.\\


\noindent\textbf {MR(2010)Subject Classification} 16T05

\section{Introduction}

Hom-algebras and Hom-coalgebras were introduced by Makhlouf and Silvestrov in
\cite{MS08} as generalizations of ordinary algebras and coalgebras in the following sense: the
associativity of the multiplication is replaced by the Hom-associativity and similar for
Hom-coassociativity. They also defined the structures of Hom-bialgebras and Hom-Hopf
algebras, and described some of their properties extending properties of ordinary bialgebras
and Hopf algebras in \cite{MS09} and \cite{MS10}.
Recently, many more properties and structures of Hom-Hopf algebras have been developed,
see \cite{CWZ13}, \cite{FG10}, \cite{GC14}, \cite{Y09}  and references cited therein.
\smallskip

Caenepeel and Goyvaerts \cite{CG11} studied Hom-bialgebras and Hom-Hopf algebras from
a categorical view point, and called them monoidal Hom-bialgebras and monoidal Hom-
Hopf algebras respectively, which are slightly different from the above Hom-bialgebras
and Hom-Hopf algebras.
In \cite{CZ13}, Chen, Wang and Zhang introduced  integrals of monoidal Hom-Hopf algebras  and  investigated the existence and
uniqueness of integrals for finite-dimensional monoidal Hom-Hopf algebras,
they also applied the integrals  to the Maschke type theorem for monoidal
Hom-Hopf algebras controlling the semisimplicity and separability of monoidal
Hom-Hopf algebras.
\smallskip

In \cite{GC14}, the first author and Chen introduced the notion of relative Hom-Hopf modules and proved that the functor $F$
from  the category of relative Hom-Hopf modules to the category of right $(A, \b)$-Hom-modules has a right
 adjoint. We also introduced the notion of total integrals for
Hom-comodule algebras which has effective applications in
investigating properties of relative Hom-Hopf modules,
and proved a Maschke type theorem for the category of relative Hom-Hopf modules as an application our theory.
\smallskip

Total integral is a vital notion in representation theories.
Doi \cite{D85} introduced the notion of total integral for any $H$-comodule algebra $A$ over a Hopf algebra $H$, which
has strong ties both to $\mathcal{M} ^{H}$ (i. e., the corepresentation of $H$) and to the representation of the
pair $(H,A)$ (i. e.,  the category of relative Hopf modules $\mathcal{M} ^{H} _{A}$),
and presented the famous necessary and sufficient criterion for the existence of a total integral.
Later, Menini and Militaru \cite{MM02} interpreted
the criterion for the existence of a total integral with the help of forgetful functors
$F: \mathcal{M} ^{H} _{A}\rightarrow \mathcal{M} _{A} $ and $G: \mathcal{M} ^{H} _{A}\rightarrow \mathcal{M} ^{H}$. Furthermore, the affineness criterion for affine
algebraic group schemes was studied by Cline et al. \cite{CPS77} and Oberst \cite{O77}.
Inspired by the ideas from Doi, Menini and Militaru, we want to investigate the
criterion for the existence of a total integral in the setting of monoidal Hom-Hopf algebras.
\smallskip

The paper is organized as follows. In Section 2, we recall some definitions and
properties relative to separable functors and monodial Hom-Hopf algebras which are needed later.
In Section 3, we present two pairs of important adjoint functors (see Proposition 3.1).
In Section 4, we show a necessary and
sufficient criterion for the existence of an integral on a monodial Hom-Hopf algebra (see Theorem 4.3).
In Section 5, we give the affineness criterion for relative Hom-Hopf
modules (see Theorem 5.7).

\section{Preliminaries}
In this section we recall some basic definitions and results related to our paper.
Throughout the paper, all algebraic systems are supposed to be over a  commutative ring $k$.
The reader is referred to  \cite{CG11} and \cite{GC14}
as general references about Hom-structures.
\smallskip

Let $\mathcal{C}$ be a category, there is a new category
$\widetilde{\mathscr{H}}(\mathcal{C})$ as follows: objects are
couples $(M, \mu)$, with $M \in \mathcal{C}$ and $\mu \in
Aut_{\mathcal{C}}(M)$. A morphism $f: (M, \mu)\rightarrow (N, \nu )$
is a morphism $f : M\rightarrow N$ in $\mathcal{C}$ such that $\nu
\circ f= f \circ \mu$.
\smallskip

Let $\mathscr{M}_k$ denotes the category of $k$-modules.
~$\mathscr{H}(\mathscr{M}_k)$ will be called the Hom-category
associated to $\mathscr{M}_k$. If $(M,\mu) \in \mathscr{M}_k$, then
$\mu: M\rightarrow M$ is obviously a morphism in
~$\mathscr{H}(\mathscr{M}_k)$. It is easy to show that
~$\widetilde{\mathscr{H}}(\mathscr{M}_k)$ =
(~$\mathscr{H}(\mathscr{M}_k),~\otimes,~(I, I),~\widetilde{a},
~\widetilde{l},~\widetilde{r}))$ is a monoidal category by
Proposition 1.1 in \cite{CG11}: the tensor product of $(M,\mu)$ and $(N,
\nu)$ in ~$\widetilde{\mathscr{H}}(\mathscr{M}_k)$ is given by the
formula $(M, \mu)\otimes (N, \nu) = (M\otimes N, \mu \otimes \nu)$.
\smallskip

Assume that $(M, \mu), (N, \nu), (P,\pi)\in
\widetilde{\mathscr{H}}(\mathscr{M}_k)$. The associativity and unit
constraints are given by the formulas
\begin{eqnarray*}
\widetilde{a}_{M,N,P}((m\o n)\o p)=\mu(m)\o (n\o \pi^{-1}(p)),\\
\widetilde{l}_{M}(x\o m)=\widetilde{r}_{M}(m\o x)=x\mu(m).
\end{eqnarray*}
An algebra in $\widetilde{\mathscr{H}}(\mathscr{M}_k)$ will be called a monoidal Hom-algebra.

{\bf Definition 2.1.}
A monoidal Hom-algebra  is an object  $(A, \alpha)\in
\widetilde{\mathscr{H}}(\mathscr{M}_k)$ together with a $k$-linear
map $m_A: A\o A\rightarrow A$ and an element $1_A\in A$ such that
\begin{eqnarray*}
&&\alpha(ab)=\alpha(a)\alpha(b); ~\a(1_A)=1_A,\\
&&\alpha(a)(bc)=(ab)\alpha(c); ~a1_A=1_Aa=\a(a),
\end{eqnarray*}
for all $a,b,c\in A$. Here we use the notation $m_A(a\o b)=ab$.
\smallskip

{\bf Definition 2.2.} A monoidal Hom-coalgebra  is
an object $(C,\g)\in \widetilde{\mathscr{H}}(\mathscr{M}_k)$
together with $k$-linear maps $\Delta:C\rightarrow C\otimes C,~~~
\D(c)=c_{(1)}\o c_{(2)}$ (summation implicitly understood) and
$\gamma:C\rightarrow C$ such that
\begin{eqnarray*}
\D(\g(c))=\g(c_{(1)})\o \g(c_{(2)}); ~~~\varepsilon(\g(c))=\varepsilon(c),
\end{eqnarray*}
and
\begin{eqnarray*}
&&\g^{-1}(c_{(1)})\otimes c_{(2)(1)}\otimes c_{(2)(2)}=c_{(1)(1)}\otimes c_{(1)(2)}\otimes\g(c_{(2)}),\\
&&\varepsilon(c_{(1)})c_{(2)}=\varepsilon(c_{(2)})c_{(1)}=\g^{-1}(c),
\end{eqnarray*}
for all $c\in C$.

{\bf Definition 2.3.} A monoidal Hom-bialgebra
$H=(H,\a,m,\eta, \Delta,\varepsilon)$  is a bialgebra in
the  category $\widetilde{\mathscr{H}}(\mathscr{M}_k)$. This means that $(H, \a,
m,\eta)$ is a monoidal Hom-algebra, $(H,\Delta,\alpha)$ is a monoidal Hom-coalgebra
and that $\D$ and $\v$ are morphisms of monoidal Hom-algebras, that is,
\begin{eqnarray*}
&&\Delta(ab)=a_{(1)}b_{(1)}\otimes a_{(2)}b_{(2)}; \Delta(1_H)=1_H\otimes 1_H,\\
&&\varepsilon(ab)=\varepsilon(a)\varepsilon(b),~\varepsilon(1_H)=1_H.
\end{eqnarray*}

{\bf Definition 2.4.} A monoidal Hom-Hopf algebra is
a monoidal Hom-bialgebra $(H, \alpha)$ together with a linear map
$S:H\rightarrow H$ in $\widetilde{\mathscr{H}}(\mathscr{M}_k)$ such
that
$$S\ast I=I\ast S=\eta\varepsilon,~S\alpha=\alpha S.$$

{\bf Definition 2.5.} Let $(A,\alpha)$ be a monoidal
Hom-algebra.  A right $(A,\alpha)$-Hom-module is an object
$(M,\mu)\in \widetilde{\mathscr{H}}(\mathscr{M}_k) $ consists of a
$k$-module and a linear map $\mu:M\rightarrow M$ together with a
morphism $\psi:M\otimes A\rightarrow M, \psi(m\c a)=m\c a$ in
$\widetilde{\mathscr{H}}(\mathscr{M}_k) $  such that
\begin{eqnarray*}
(m\c a)\c\alpha(b)=\mu(m)\c (ab);~~~~ m\c 1_A=\mu(m),
\end{eqnarray*}
for all $a\in A$ and $m\in M$. The fact that $\psi\in \widetilde{\mathscr{H}}(\mathscr{M}_k)$ means
\begin{eqnarray*}
\mu(m\c a)=\mu(m)\c\alpha(a).
\end{eqnarray*}
A morphism $f: (M, \mu)\rightarrow (N, \nu)$ in $\widetilde{\mathscr{H}}(\mathscr{M}_k)$ is called right $A$-linear if it preserves the $A$-action, that is, $f(m\c a)=f(m)\c a$. $\widetilde{\mathscr{H}}(\mathscr{M}_k)_A$ will denote the category of right $(A, \a)$-Hom-modules and $A$-linear morphisms.

{\bf Definition 2.6.} Let $(C,\g)$ be a monoidal
Hom-coalgebra. A right $(C,\g)$-Hom-comodule  is an object
$(M,\mu)\in \widetilde{\mathscr{H}}(\mathscr{M}_k)$  together with a
$k$-linear map $\rho_M: M\rightarrow M\otimes C$
($\rho_M(m)=m_{[0]}\o m_{[1]}$) in
$\widetilde{\mathscr{H}}(\mathscr{M}_k)$ such that
\begin{eqnarray*}
m_{[0][0]}\otimes (m_{[0][1]}\otimes
\g^{-1}(m_{[1]}))=\mu^{-1}(m_{[0]})\otimes \D_C(m_{[1]});
~m_{[0]}\varepsilon(m_{[1]})=\mu^{-1}(m),
\end{eqnarray*}
for all $m\in M$.  The fact that $\rho_M\in \widetilde{\mathscr{H}}(\mathscr{M}_k)$ means
\begin{eqnarray*}
\rho_M(\mu(m))=\mu(m_{[0]})\otimes \g(m_{[1]}).
\end{eqnarray*}
Morphisms of right $(C, \g)$-Hom-comodule are defined in the obvious way. The
category of right $(C, \g)$-Hom-comodules will be denoted by $\widetilde{\mathscr{H}}(\mathscr{M}_k)^C$.

{\bf Definition 2.7.}  Let $(H,\alpha)$  be a monoidal Hom-Hopf algebra.  A monoidal Hom-algebra  $(A,\beta)$ is called  a right
$(H,\alpha)$-Hom-comodule algebra if $(A,\beta)$ is a right $(H,\alpha)$ Hom-comodule with a coaction $\rho_A: A\rightarrow A\o H,
 \rho_A(a)= a_{[0]}\o a_{[1]}$
 such that the following conditions hold:
 \begin{eqnarray*}
 &&\rho_{A} ( ab)= a_{[0]} b_{[0]} \otimes a_{[1]}b_{[1]},
 \rho_{A} ( 1_A)=1_A\o 1_H,
 \end{eqnarray*}
for all $a,b \in A$.

\section{Relative Hom-Hopf Modules: functors and structures}
\def\theequation{\arabic{section}.\arabic{equation}}
\setcounter{equation} {0}
Let $(H,\alpha)$  be a monoidal Hom-Hopf algebra and  $(A,\b)$  a right
$(H,\a)$-Hom-comodule algebra.
Recall from \cite{GC14},
a relative Hom-Hopf module $(M,\mu)$
is a right $(A,\b)$-Hom-module which is also a right
$(H,\a)$-Hom-comodule with the  coaction structure $ \rho _{M} : M
\rightarrow M \otimes  H $ defined by $\rho _{M}(m)=m_{[0]}\o
m_{[1]}$ such that the following compatible condition holds: for all
$m\in M $ and $a\in A$,
\begin{eqnarray*}
&&\rho _{M}( m a) = m _{[0]} \c a _{[0]} \otimes  m _{[1]} a _{[1]}.
\end{eqnarray*}

A morphism between two right relative Hom-Hopf modules is a $k$-linear map
which is a morphism in the categories $\widetilde{\mathscr{H}}(\mathscr{M}_k)_A$ and $\widetilde{\mathscr{H}}(\mathscr{M}_k)^H$ at the same time.  $ \widetilde{\mathscr{H}}(\mathscr{M}_k)^{H}_{A}$ will denote  the category of right relative Hom-Hopf modules and morphisms
between them.

By \cite{GC14}, the forgetful functor
$F:  \widetilde{\mathscr{H}}(\mathscr{M}_k)^{H}_{A} \rightarrow \widetilde{\mathscr{H}}(\mathscr{M}_k)_{A}$ has a right
 adjoint $G:  \widetilde{\mathscr{H}}(\mathscr{M}_k)_{A}\rightarrow \widetilde{\mathscr{H}}(\mathscr{M}_k)^{H}_{A}$.
 $G$ defined by
$$
 G(M) = M \otimes H
 $$
 with structure maps
\begin{eqnarray*}
&&( m  \otimes h  ) \cdot a =   m \c a _{[0]}\otimes h  a _{[1]},\\
&&\rho _{G(M)} ( m  \otimes h ) =  (\mu^{-1}( m ) \otimes h _{(1)}) \otimes \alpha(h _{(2)}),
\end{eqnarray*}
for all $a \in A, m\in M$ and $h\in H$.

Let us describe the
 unit $\eta$ and the counit $\delta$ of the adjunction. The unit is described by the coaction:
 for $M \in   \widetilde{\mathscr{H}}(\mathscr{M}_k)^{H}_{A} $,  we define
 $\eta _{M} :
  M \rightarrow M \otimes H$  as follows:
 $$
 \eta _{M}(m) =  m _{[0]}
 \otimes m _{[1]},~m \in M.
 $$
We can check that $\eta_{M} \in  \widetilde{\mathscr{H}}(\mathscr{M}_k)^{H}_{A}$ .
For any $N \in \widetilde{\mathscr{H}}(\mathscr{M}_k)_{A}$, define
 $$
 \delta _{N}: N\otimes H\rightarrow N,~\delta _{N}(n  \otimes h  )
 = \varepsilon(h )n, ~n \in N,h \in H .
 $$
It is easy to check that $\eta _{M}\in   \widetilde{\mathscr{H}}(\mathscr{M}_k)^{H}_{A}$ .
One may verify that $\eta$ and $\delta$ defined above are natural transformations and
 satisfy
 $$
 G (\delta _{N})\circ \eta _{G (N)} = I _{G (N)},~
 \delta _{F (M)}\circ F (\eta _{M}) = I _{F (M)},
 $$
 for all $M \in  \widetilde{\mathscr{H}}(\mathscr{M}_k)^{H}_{A}$ and $N \in  \widetilde{\mathscr{H}}(\mathscr{M}_k)_{A}$ .

Since $(A,\b)$ is a right $(A,\b)$-Hom-module, $G(A)$ is a relative Hom-Hopf module via

\begin{eqnarray}
&&(a \otimes h ) \cdot b = a b _{[0]}\otimes h b _{[1]};\\
&&\rho_{A\o H}  ( a \otimes h  )
  =(\b^{-1}(a) \otimes h _{(1)})\otimes \a(h _{(2)}),
\end{eqnarray}
for all $a, b\in A$ and $h\in H$. Notice that
$
GFG(A)=A\o H\o H
$
is an object in $\widetilde{\mathscr{H}}(\mathscr{M}_k)^H$, and we also have
$$
\left\{
  \begin{array}{ll}
(a \otimes h \o g ) \cdot b = a b _{[0][0]}\otimes h b _{[0][1]}\o gb_{[1]};\\
  \rho_{A\o H\o H}  ( a \otimes h \o g )
  =(\b^{-1}(a) \o \a^{-1}(h)\otimes g _{(1)})\otimes \a(g _{(2)}),
\end{array}
\right.
$$
for any $a, b\in A$ and $h, g\in H$.

Similar to \cite{GC14}, one can  check that $\widetilde{F}:  \widetilde{\mathscr{H}}(\mathscr{M}_k)^{H}_{A} \rightarrow \widetilde{\mathscr{H}}(\mathscr{M}_k)^{H}$
which forgets the $(A,\b)$-action has a left adjoint given by
\begin{eqnarray*}
\widetilde{G}:  \widetilde{\mathscr{H}}(\mathscr{M}_k)^{H}\rightarrow \widetilde{\mathscr{H}}(\mathscr{M}_k)^{H}_{A}, \widetilde{G}(N)=A\o N,~N\in \widetilde{\mathscr{H}}(\mathscr{M}_k)^{H}.
\end{eqnarray*}
The relative Hom-Hopf module structure on $\widetilde{G}(N)$ is given by
$$
\left\{
  \begin{array}{ll}
    (a \otimes n) \cdot b = a\b(b) \otimes \nu^{-1}(n), \\
  \rho_{A\o N} ( a \otimes n)= (a _{[0]}\o n _{[0]} \otimes )\otimes a _{[1]}n _{[1]},
  \end{array}
\right.
$$
for all $a\in A$ and $n\in N$. Since $(H,\a)$ is a right
$(H,\a)$-Hom-comodule, then $\widetilde{G}(H)=A\o H$ can be viewed as a relative Hom-Hopf
module via
\begin{eqnarray}
&&(a \otimes h) \cdot b = a\b(b) \otimes \a^{-1}(h), \\
&&\rho_{A\o H} ( a \otimes h)= (a _{[0]}\o h _{(1)} \otimes )\otimes a _{[1]}h _{(2)},
\end{eqnarray}
for all $a\in A$ and $h\in H$.
\smallskip

From Eqs. (3.1) and (3.2) or Eqs. (3.3) and (3.4), we have two types of relative Hom-Hopf modules structures on $A\o H$. The following proposition reveals their
relationship.

{\bf Proposition 3.1.} With the two types of relative Hom-Hopf modules on $A\o H$.
Then
\begin{eqnarray*}
G(A)\cong\widetilde{G}(H)
\end{eqnarray*}
 in $\widetilde{\mathscr{H}}(\mathscr{M}_k)^{H}_{A}$.

{\bf Proof.} For any $a\in A$ and $h\in H$, we construct two maps $u$  and $v$  as follows:
\begin{eqnarray*}
&&u: A\o H\rightarrow A\o H,  u(a\o h)=\b(a_{[0]})\o h\a^{-1}(a_{[1]}),\\
&&v: A\o H\rightarrow A\o H,  u(a\o h)=\b(a_{[0]})\o h\a^{-1}(S(a_{[1]})).
\end{eqnarray*}
We claim that $u\circ v=id_{A\o H}$. In fact,
\begin{eqnarray*}
u\circ v(a\o h)&=&u(\b(a_{[0]})\o h\a^{-1}(S(a_{[1]})))\\
&=&\b^2(a_{[0][0]})\o (h\a^{-1}(S(a_{[1]})))\a^{-1}(a_{[0][1]})\\
&=&\b(a_{[0]})\o (hS(a_{[1](2)}))\a^{-1}(a_{[1](1)})\\
&=&\b(a_{[0]})\o \a^{-1}(h)(S(a_{[1](2)})a_{[1](1)})\\
&=&\b(a_{[0]})\o \a^{-1}(h)\varepsilon(a_{[1]})\\
&=&a\o h.
\end{eqnarray*}
Similarly, one may verify that $v\circ u=id_{A\o H}$.
 The proof that $u$ is a
morphism in $\widetilde{\mathscr{H}}(\mathscr{M}_k)^{H}_{A}$ is straightforward.$\hfill \Box$

\section{Natural transformations versus integrals}
\def\theequation{4. \arabic{equation}}
\setcounter{equation} {0} \hskip\parindent

{\bf Definition 4.1.}$^{\cite{GC14}}$ Let $(H, \a)$ be a monoidal Hom Hopf algebra and  $(A,\b)$  a right
$(H,\a)$-Hom-comodule algebra. The map $\phi: (H, \a)\rightarrow (A,\b)$ is called a  total integral if the following conditions
are satisfied:
$$\rho_A\phi=(\phi\o id_H)\Delta_H,~~~ \phi\a=\b\phi,~~~~~\phi(1_H)=1_A.$$

{\bf Proposition 4.2.}  Let $(A,\b)$ be a right
$(H,\a)$-Hom-comodule algebra,
 $(M, \mu), (N, \nu) \in \widetilde{\mathscr{H}}(\mathscr{M}_k)^{H}_{A}$
 and $u: (N, \nu) \rightarrow (M, \mu) $ a $k$-linear map. Suppose that there exists a
total integral $\phi: (H, \a)\rightarrow (A,\b)$. Then

(1) The map \begin{eqnarray*}
           \widetilde{u}: (N,\nu)\rightarrow (M,\mu),~\widetilde{u}(n)=u(\nu(n_{[0]}))_{[0]}\c \phi(S(u(n_{[0]})_{[1]})\a^{-1}(n_{[1]})),~n\in N
            \end{eqnarray*}
is right $(H, \a)$-colinear.

(2) If $f : (M, \mu) \rightarrow (N, \nu)$ is a morphism in $\widetilde{\mathscr{H}}(\mathscr{M}_k)^{H}_{A}$ of which is a $k$-split
injection (respectively, a $k$-split surjection), then $f$ has an $(H, \a)$-colinear retraction
(respectively, a section).

{\bf Proof. } For any $n\in N$, we have
\begin{eqnarray*}
\rho_M\circ \widetilde{u}(n)&=&u(\nu(n_{[0]}))_{[0][0]}\c \phi(S(u(n_{[0]})_{[1]})\a^{-1}(n_{[1]}))_{[0]}\\
&&\hspace{5cm}\o u(\nu(n_{[0]}))_{[0][1]} \phi(S(u(n_{[0]})_{[1]})\a^{-1}(n_{[1]}))_{[1]}\\
&=&u(n_{[0]})_{[0]}\c \phi(S(u(\nu(n_{[0]}))_{[1](2)})\a^{-1}(n_{[1]}))_{[0]}\\
&&\hspace{4cm}\o u(\nu(n_{[0]}))_{[1](1)} \phi(S(u(\nu(n_{[0]}))_{[1](2)})\a^{-1}(n_{[1]}))_{[1]}\\
&=&u(n_{[0]})_{[0]}\c \phi(S(u(\nu(n_{[0]}))_{[1](2)(2)})\a^{-1}(n_{[1](1)}))\\
&&\hspace{4cm}\o u(\nu(n_{[0]}))_{[1](1)}S(u(\nu(n_{[0]}))_{[1](2)(1)})\a^{-1}(n_{[1](2)})\\
&=&u(n_{[0]})_{[0]}\c \phi(S(u(n_{[0]})_{[1](2)})\a^{-1}(n_{[1](1)}))\\
&&\hspace{4cm}\o u(\nu^2(n_{[0]}))_{[1](1)(1)} S(u(\nu(n_{[0]}))_{[1](1)(2)})\a^{-1}(n_{[1](2)})\\
&=&u(n_{[0]})_{[0]}\c \phi(S(u(n_{[0]})_{[1](2)})\a^{-1}(n_{[1](1)}))\\
&&\hspace{4cm}\o [u(\nu(n_{[0]}))_{[1](1)(1)} S(u(\nu(n_{[0]}))_{[1](1)(2)})]n_{[1](2)}\\
&=&u(n_{[0]})_{[0]}\c \phi(S(u(\nu^{-1}(n_{[0]}))_{[1]})\a^{-1}(n_{[1](1)}))\o \a(n_{[1](2)})\\
&=&u(\nu(n_{[0][0]}))_{[0]}\c \phi(S(u(n_{[0][0]})_{[1]})\a^{-1}(n_{[1](1)}))\o n_{[1](2)}\\
&=&( \widetilde{u}\o id_H)\circ \rho_N(n).
\end{eqnarray*}
Hence $\widetilde{u}$ is  right $(H, \a)$-colinear.

Let $u: (N, \nu)\rightarrow (M, \mu)$ be $k$-linear, then $\widetilde{u}: (N, \nu)\rightarrow (M, \mu)$ is  right
$(H, \a)$-colinear. Assume that $u$ is a retraction of $f$, for any $m\in M$, we have
\begin{eqnarray*}
(\widetilde{u}\circ f)(m)&=& u(\nu(f(m)_{[0]}))_{[0]}\c \phi(S(u(f(m)_{[0]})_{[1]})\a^{-1}(f(m)_{[1]}))\\
&=& u(\nu(f(m_{[0])}))_{[0]}\c \phi(S(u(f(m_{[0]}))_{[1]})\a^{-1}(m_{[1]}))\\
&=& \nu(m_{[0][0]})\c \phi(S(m_{[0][1]})\a^{-1}(m_{[1]}))\\
&=& m_{[0]}\c \phi(S(m_{[1](1)})m_{[1](2)})\\
&=&m_{[0]}\varepsilon(m_{[1]})\c 1_A=m.
\end{eqnarray*}
So $\widetilde{u}: (N, \nu)\rightarrow (M, \mu)$ is a right
$(H, \a)$-colinear retraction of $f$. On the other hand, if $u$
is a section of $f$, for any $n\in N$, we have
\begin{eqnarray*}
(f \circ\widetilde{u})(n)&=& f(u(\nu(n_{[0]}))_{[0]})\c \phi(S(u(n_{[0]})_{[1]})\a^{-1}(n_{[1]}))\\
&=& f(u(\nu(n_{[0]}))_{[0]})\c \phi(S(u(n_{[0]})_{[1]})\a^{-1}(n_{[1]}))\\
&=& \nu(n_{[0][0]}) \c \phi(S(n_{[0][1]})\a^{-1}(n_{[1]}))\\
&=& n_{[0]} \c \phi(S(n_{[1](1)})n_{[1](2)})\\
&=&n_{[0]}\varepsilon(n_{[1]})\c 1_A=n.
\end{eqnarray*}
Thus $\widetilde{u}: (N, \nu)\rightarrow (M, \mu)$ is a right $H$-colinear section of $f$.$\hfill \Box$
\smallskip

{\bf Theorem 4.3.} Let $(A,\b)$ be a right
$(H,\a)$-Hom-comodule algebra. Then the following statements
are equivalent:

(1) There exists a total integral $\phi: (H, \a)\rightarrow (A,\b)$;

(2) There exists a natural transformation  $\rho: \widetilde{F}\circ 1_{\widetilde{\mathscr{H}}(\mathscr{M}_k)^{H}_{A}}\rightarrow \widetilde{F}\circ G\circ F$ splits;

(3) $\rho_A: A\rightarrow A\o H$ splits in $\widetilde{\mathscr{H}}(\mathscr{M}_k)^{H}$.

 Consequently, if one
of the equivalent conditions holds, then any relative Hopf module is injective
as a right $(H,\a)$-comodule.

{\bf Proof.} (1)$\Rightarrow$(2) Assume that $\phi: (H, \a)\rightarrow (A,\b)$ is a total integral.
We have to
construct a natural transformation $\lambda$ that splits $\rho$.
Let $(M, \mu)\in \widetilde{\mathscr{H}}(\mathscr{M}_k)^{H}_{A}$ and $\varphi: M\o H\rightarrow M$
be the $k$-linear retraction of $\rho_M: M\rightarrow M\o H,~\varphi(m\o h)=\varepsilon(h)m$.
Define
 \begin{eqnarray*}
             &&\lambda_M: \widetilde{F}\circ G\circ F(M)\rightarrow  \widetilde{F}(M),\\
             &&\lambda_M(m\o h)=\mu(m_{[0]})\c \phi(S(m_{[1]})\a^{-1}(h)).
              \end{eqnarray*}
It follows from Proposition 4.2 that the map $\lambda_M$ is a right $(H, \a)$-colinear
retraction of $\rho_M$.

It remains to prove that $\lambda_M\in \widetilde{\mathscr{H}}(\mathscr{M}_k)^{H}_{A}$ is a natural transformation.
Let $f: (M,\mu)\rightarrow (N,\nu)$ be a morphism in $\widetilde{\mathscr{H}}(\mathscr{M}_k)^{H}_{A}$, we have to prove that
\begin{eqnarray*}
f\circ \lambda_M=\lambda_N\circ(f\o id_H).
\end{eqnarray*}
Since $f$ is right $(H, \a)$-colinear, we have                                                                                                                \begin{eqnarray*}
\lambda_N\circ(f\o id_H)(m\o h)&=&\lambda_N(f(m)\o h)\\
&=&\nu(f(m)_{[0]})\c \phi(S(f(m)_{[1]})\a^{-1}(h))\\
&=&\nu(f(m_{[0]}))\c \phi(S(m_{[1]})\a^{-1}(h))\\
&=&f\circ \lambda_M(m\o h).
\end{eqnarray*}
So $\lambda_M$ is a natural transformation that splits $\rho_M$.
\smallskip

(2)$\Rightarrow$(3)  Assume that for any $(M, \mu)\in \widetilde{\mathscr{H}}(\mathscr{M}_k)^{H}_{A}$, $\chi_M: \widetilde{F}(M)\rightarrow \widetilde{F}\circ G\circ F(M)$ splits in the category $\widetilde{\mathscr{H}}(\mathscr{M}_k)^{H}$ of right $H$-comodules and the character
of the splitting is functorial. In particular, we consider the relative Hopf module $(A,\b)$, then we have
\begin{eqnarray*}
\rho_A: A\rightarrow A\o H, ~~~~~\rho_A(a)=a_{[0]}\o a_{[1]}
\end{eqnarray*}
splits in $\widetilde{\mathscr{H}}(\mathscr{M}_k)^{H}$ and let
\begin{eqnarray*}
\lambda_A: A\o H\rightarrow A
\end{eqnarray*}
be a right $(H, \a)$-colinear retraction of it.

(3)$\Rightarrow$(1) Assume that the map $\rho_A: A\rightarrow A\o H,  \rho_A(a)=a_{[0]}\o a_{[1]}$ splits in $\widetilde{\mathscr{H}}(\mathscr{M}_k)^{H}$.
Let $\lambda_A: A\o H\rightarrow A$  be the map which splits $\rho_A$. We define
\begin{eqnarray*}
\phi: (H, \a)\rightarrow (A,\b), ~~~~\phi(h)=i_A^{-1}(\lambda_A(1_A\o h)),
\end{eqnarray*}
where $i_A: A\rightarrow A$ is the natural isomorphism.
Since $\lambda_A$ is right $(H, \a)$-colinear, we have
\begin{eqnarray*}
\lambda_A(1_A\o h_{(1)})\o h_{(2)}=\lambda_A(1_A\o h)_{[0]}\o \lambda_A(1_A\o h)_{[1]}.
\end{eqnarray*}
Let $\lambda_A(1_A\o h_{(1)})=b$ and $\lambda_A(1_A\o h)=q$,  we get
$
b\o h_{(2)}=q_{[0]}\o q_{[1]}.
$
Therefore,  $\phi$ is $(H, \a)$-colinear. It is sufficient to
check that
$
\rho_A\phi=(\phi\o id_H)\Delta_H.
$
For any  $h\in H$, we have
\begin{eqnarray*}
\rho_A\circ \phi(h)&=&\rho_A(i_A^{-1}(\lambda_A(1_A\o h)))
=q_{[0]}\o q_{[1]}.
\end{eqnarray*}
On the other hand,
\begin{eqnarray*}
&&(\phi\o id_H)\Delta_H(h)=\phi(h_{(1)})\o h_{(2)}\\
&=&i_A^{-1}(\lambda_A(1_A\o h_{(1)}))\o h_{(2)}
= b\o h_{(2)}.
\end{eqnarray*}
Notice that
$
\phi(1_H)=i_A^{-1}(\lambda_A(1_A\o 1_H))
=i_A^{-1}(\lambda_A(\rho_A(1_A)))=i_A^{-1}(1_A)=1_A.
$
It is easy to check that $\phi\a=\b\phi$. So $\phi$ is a  total integral.
And this completes the proof.
$\hfill \Box$
\smallskip

From Theorem 4.3, we have the following corollary.
\smallskip

{\bf Corollary 4.4.} Let $F:  \widetilde{\mathscr{H}}(\mathscr{M}_k)^{H}_{A} \rightarrow \widetilde{\mathscr{H}}(\mathscr{M}_k)_{A}$ and $\widetilde{F}:  \widetilde{\mathscr{H}}(\mathscr{M}_k)^{H}_{A} \rightarrow \widetilde{\mathscr{H}}(\mathscr{M}_k)^{H}$
be two forgetful functors.  Then $F$ is $\widetilde{F}$-separable if and
only if there exists a total integral $\phi: (H, \a)\rightarrow (A,\b)$.
\smallskip

{\bf Definition 4.5.}  Let $(H, \a)$ be a monoidal Hom-Hopf algebra with a bijective antipode and $(A,\b)$  a right
$(H,\a)$-Hom-comodule algebra. A $k$-linear map $\gamma: H\rightarrow Hom(H,A)$
satisfying $\gamma(\a(g))(\a(h))=\b\circ \gamma(g)(h)$ is called a quantum integral if
\begin{eqnarray}
 && \gamma(\a^{-1}(g))( h _{(1)})\otimes \a(h _{(2)}) = [\gamma( \a(g _{(2)}))(h)]_{[0]}\otimes g _{(1)} [\gamma(g _{(2)}(\a^{-1}(h))]_{[1]},
 \end{eqnarray}
for all $g,h\in H$. A quantum integral $\gamma: H\rightarrow Hom(H,A)$ is called total if
\begin{eqnarray}
\gamma(h_{(1)})(h_{(2)})=\varepsilon(h)1_A,
\end{eqnarray}
for all $h\in H$.
\smallskip

{\bf Remark 4.6.} (1) Let $\gamma: H\rightarrow Hom(H,A)$ be a quantum integral, then
\begin{eqnarray*}
\phi: (H, \a)\rightarrow (A,\b), ~\phi(h)=\gamma(h)(1_H)
\end{eqnarray*}
satisfies the condition:
\begin{eqnarray*}
\phi(h_{(1) })\otimes h _{(2)}=\phi(h)_{[0]} \otimes \phi  (h)_{[1]},
\end{eqnarray*}
for all $h\in H$, that is, $\phi: (H, \a)\rightarrow (A,\b)$ is right $(H,\a)$-colinear.
\smallskip

(2) As opposed to the case $\widetilde{\mathscr{H}}(\mathscr{M}_k)^{H}_{A}$, a right $(H,\a)$-colinear map $\phi: (H, \a)\rightarrow (A,\b)$ is not sufficient to construct a quantum integral. If $\phi: (H, \a)\rightarrow (A,\b)$ is a $k$-linear map satisfying a more powerful condition:
\begin{eqnarray*}
g \phi(h)_{[1]} \otimes \phi(h)_{[0]} = \phi(h)_{[1]}g \otimes \phi(h)_{[0]},~h, g\in H.
\end{eqnarray*}
Then
$
\gamma: H\rightarrow Hom(H,A), ~\gamma(g)(h)=\phi(hS^{-1}(g))
$
is a quantum integral.
\smallskip

{\bf Example 4.7.} (1) Recall from that let $H=M^n(k)$ be the $n\times n$ matrix monoidal Hom-Hopf algebra with the structures given by
\begin{eqnarray*}
c_{iu}c_{uj}=\delta_{ij}\a(c_{ij}), ~~~~~\sum c_{ii}=1_H,\\
\Delta(c_{ij})=\a^{-1}(c_{iu})\o \a^{-1}(c_{uj}),~~~~~\varepsilon(c_{ij})=\delta_{ij}.
\end{eqnarray*}
Let $(M^n(k),A)$ be a relative Hom-Hopf datum and $B=A^{coH}=\{a\in A|\rho(a)=\b^{-1}(a)\o 1_H\}$ the subalgebra of coinvariants of $A$.
 Let $\mu=(\mu_{ij})\in M^n(B)$ be an arbitrary $n\times n$-matrix over $B$.
 Then the map
\begin{eqnarray*}
\gamma: H\rightarrow Hom(H,A),~\g(c_{ij})(c_{rs})=\delta_{is}\b(\mu_{rj})
\end{eqnarray*}
is an integral of $(M^n(k),A)$.
Indeed, for any $h=c_{ij}$ and $g=c_{kl}$, we have
\begin{eqnarray*}
 \gamma(\a^{-1}(g))( h _{(1)})\otimes \a(h _{(2)})= \gamma(\a^{-1}(c_{kl}))( \a^{-1}(c_{iu}))\otimes c_{uj}=\delta_{ku}\mu_{il}\otimes c_{uj}=\mu_{il}\otimes c_{kj}
\end{eqnarray*}
and
\begin{eqnarray*}
&&[\gamma( \a(g _{(2)}))(h)]_{[0]}\otimes g _{(1)} [\gamma(g _{(2)}(\a^{-1}(h))]_{[1]}\\
&=&[\gamma( c_{ul})(c_{ij})]_{[0]}\otimes \a^{-1}(c_{ku}) [\gamma(\a^{-1}(c_{ul}))(\a^{-1}(c_{ij}))]_{[1]}\\
&=&\delta_{uj}\mu_{il}\otimes c_{ku}=\mu_{il}\otimes c_{kj}.
\end{eqnarray*}
So $\g$ is an integral. On the other hand, we have
\begin{eqnarray*}
\gamma(h_{(1)})(h_{(2)})=\gamma(c_{iu})(c_{uj})=\delta_{ij}\b(\mu_{uu}).
\end{eqnarray*}
Hence $\g$ is a total integral if and only $\b(\mu_{uu})=1_A$.
\smallskip

(2) Another class of examples of total integrals arises from the
graded case. Recall from \cite{CZ13} that let $G$ be a finite group and $H=kG$ be the monoidal Hom-Hopf algebra with structures:
\begin{eqnarray*}
&&xy=\a(xy), 1x=\a(x),\varepsilon(x)=1\\
&&\Delta(x)=\a^{-1}(x)\o \a^{-1}(x),~x,y\in G.
\end{eqnarray*}
Let $(kG,A)$ be a relative Hom-Hopf datum and $\mu=(\mu_{xy})_{x,y\in G}$ a family of elements of $B=A^{coH}=\{a\in A|\rho(a)=\b^{-1}(a)\o 1_H\}$.
Then the map
\begin{eqnarray*}
\gamma: H\rightarrow Hom(H,A),~\g(x)(y)=\delta_{xy}\b(\mu_{xy})
\end{eqnarray*}
is an integral of $(kG, A)$. Indeed, for any $h=x$ and $g=y$, we have
\begin{eqnarray*}
&&[\gamma(y)(x)]_{[0]}\otimes \a^{-1}(y) [\gamma(\a^{-1}(y)(\a^{-1}(x))]_{[1]}\\
&=& \delta_{xy}\mu_{xy}\o y=\delta_{xy}\mu_{xy})\o x=\gamma(\a^{-1}(y))( \a^{-1}(x))\otimes x.
\end{eqnarray*}
Furthermore, $\g$ is a total integral iff $\mu_{xx}=1_A$ for all $x\in G$.
\smallskip

Next we will give some applications of the existence of total quantum integrals.
\smallskip

{\bf Theorem 4.8.} Let $(A,\b)$ be a right
$(H,\a)$-Hom-comodule algebra and $(M, \mu)$ an object in $\widetilde{\mathscr{H}}(\mathscr{M}_k)^{H}_{A}$.
Assume that there exists a total quantum  integral $\gamma: H\rightarrow Hom(H,A)$,
 then the map
\begin{eqnarray*}
&&f: A\o H\o M\rightarrow M,\\
&&f(a\o h\o m)=\mu(m_{[0]})\gamma(\a^{-1}(m_{[1]}))(\a^{-2}(h)S^{-1}(\a^{-1}(a_{[1]})))\b(a_{[0]}),
\end{eqnarray*}
for all $h\in H, a\in A$ and $m\in M$, is a $k$-split epimorphism in $\widetilde{\mathscr{H}}(\mathscr{M}_k)^{H}_{A}$.
In particular, $A\o H$ is a generator in the category $\widetilde{\mathscr{H}}(\mathscr{M}_k)^{H}_{A}$.

{\bf Proof.} We first show that $A\o H\o M$ is an object in $ \widetilde{\mathscr{H}}(\mathscr{M}_k)^{H}_{A}$
with structures arising from the ones of $A\o H$, ie.,
 \begin{eqnarray*}
&&(a\o h\o m)b= a\b^{-1}(b_{[0]})\o h\a^{-1}(b_{[1]})\o \mu(m),\\
&&\rho_{A\o H\o M}(a\o h\o m)= \b^{-1}(a)\o h_{(1)}\o \mu^{-1}(m) \o \a^{2}(h_{(2)}),
 \end{eqnarray*}
for all $h\in H, a,b\in A$ and $m\in M$.
It is routine to check that $A\o H\o M$ is a  right $(H,\a)$-Hom-comodule and a right $(A,\b)$-Hom-module.
Now we only check the compatibility condition.
 Indeed, for any $a, b \in A $, we have
\begin{eqnarray*}
 &  & \rho_{A\o H\o M}  ((a  \otimes h\o  m) b)\\
 & =&  \rho_{A\o H\o M}  ( a\b^{-1}(b_{[0]})\o h\a^{-1}(b_{[1]})\o \mu(m))\\
 & = &  \b^{-1}(a)\b^{-2}(b_{[0]})\o h_{(1)}\a^{-1}(b_{[1](1)})\o m\o \a^{2}(h_{(2)})\a(b_{[1](2)})\\
 & = & \b^{-1}(a)\b^{-1}(b_{[0][0]})\o h_{(1)}\a^{-1}(b_{[0][1]})\o m\o \a^{2}(h_{(2)})b_{[1]}\\
 & = & [\b^{-1}(a)\o h_{(1)}\o \mu^{-1}(m)]b_{[0]}\o \a^{2}(h_{(2)})b_{[1]}\\
 &=& \rho_{A\o H\o M} (a  \otimes h\o  m )b,
  \end{eqnarray*}
as desired.
Next we will prove that $f$ is a $k$-split surjection.
 Let $g: M\rightarrow A\o H \o M, g(m)= 1_A\o \a^{-1}(m_{[1]})\o m_{[0]}$, for all $m\in M$.
  Then $g$ is right $(H, \a)$-colinear and for any $m\in M$, we have
\begin{eqnarray*}
(f\circ g)(m)&=& f(1_A\o \a^{-1}(m_{[1]})\o m_{[0]})\\
&=&\mu(m_{[0][0]})\gamma(m_{[0][1]})(\a^{-1}(m_{[1]}))\\
&=&m_{[0]}\gamma(m_{[1](1)})(m_{[1](2)})\\
&=&m.
\end{eqnarray*}
For all $a,b\in A$ and $h\in H, m\in M$, we have
\begin{eqnarray*}
&&f((a\o h\o m)b)\\
&=& f(a\b^{-1}(b_{[0]})\o h\a^{-1}(b_{[1]})\o \mu(m))\\
&=& \mu^{2}(m_{[0]})\gamma(m_{[1]})(\a^{-2}(h\a^{-1}(b_{[1]}))S^{-1}(\a^{-1}(a_{[1]}\a^{-1}(b_{[0][1]}))))\b(a_{[0]}\b^{-1}(b_{[0][0]}))\\
&=& \mu^{2}(m_{[0]})\gamma(m_{[1]})(\a^{-2}(h)\a^{-3}(b_{[1]})S^{-1}(\a^{-2}(b_{[0][1]}))S^{-1}(\a^{-1}(a_{[1]})))[\b(a_{[0]})b_{[0][0]}]\\
&=& \mu^{2}(m_{[0]})\gamma(m_{[1]})(\a^{-1}(h)\{[\a^{-3}(b_{[1](2)})S^{-1}(\a^{-3}(b_{[1](1)}))]S^{-1}(\a^{-1}(a_{[1]}))\})[\b(a_{[0]})\b^{-1}(b_{[0]})]\\
&=& \mu^{2}(m_{[0]})\gamma(m_{[1]})(\a^{-1}(h)S^{-1}(a_{[1]}))[\b(a_{[0]})\b^{-2}(b)]\\
&=& \mu^{2}(m_{[0]})\{[\gamma(\a^{-1}(m_{[1]}))(\a^{-2}(h)S^{-1}(\a^{-1}(a_{[1]})))\b(a_{[0]})]\b^{-1}(b)\}\\
&=& \{\mu(m_{[0]})[\gamma(\a^{-1}(m_{[1]}))(\a^{-2}(h)S^{-1}(\a^{-1}(a_{[1]})))\b(a_{[0]})]\}b\\
&=&f(a\o h\o m)b.
\end{eqnarray*}
Then $f$ is a right $(A,\b)$-linear. It remains to prove that $f$ is also right $(H,\a)$-colinear.
For this purpose, we do the following calculations:
\begin{eqnarray*}
&&( f \o id_H)\rho_{A\o H\o M}(a\o h\o m)\\
&=&( f \o id_H)( \b^{-1}(a)\o h_{(1)}\o \mu^{-1}(m) \o \a^{2}(h_{(2)}))\\
&=& f(\b^{-1}(a)\o h_{(1)}\o \mu^{-1}(m) )\o \a^{2}(h_{(2)}))\\
&=& m_{[0]}\gamma(\a^{-2}(m_{[1]}))(\a^{-2}(h_{(1)})S^{-1}(\a^{-2}(a_{[1]})))a_{[0]} \o \a^{2}(h_{(2)})).
\end{eqnarray*}
On the other side,
\begin{eqnarray*}
&&\rho_{A\o H\o M}\circ f(a\o h\o m)\\
&=& \rho_{A\o H\o M}(\mu(m_{[0]})\gamma(\a^{-1}(m_{[1]}))(\a^{-2}(h)S^{-1}(\a^{-1}(a_{[1]})))\b(a_{[0]}))\\
&=& \mu(m_{[0][0]})\gamma(\a^{-1}(m_{[1]}))(\a^{-2}(h)S^{-1}(\a^{-1}(a_{[1]})))_{[0]}\b(a_{[0][0]})\\
&&\o\mu(m_{[0][1]})\gamma(\a^{-1}(m_{[1]}))(\a^{-2}(h)S^{-1}(\a^{-1}(a_{[1]})))_{[1]}\b(a_{[0][1]})\\
&=& m_{[0]}\gamma(m_{[1](2)})(\a^{-2}(h)S^{-1}(\a^{-1}(a_{[1]})))_{[0]}\b(a_{[0][0]})\\
&&\o\mu(m_{[1](1)})\gamma(m_{[1](2)})(\a^{-2}(h)S^{-1}(\a^{-1}(a_{[1]})))_{[1]}\b(a_{[0][1]})\\
&=& m_{[0]}\gamma(m_{[1](2)})(\a^{-2}(h)S^{-1}(\a^{-1}(a_{[1]})))_{[0]}\b(a_{[0][0]})\\
&&\o [m_{[1](1)}\gamma(m_{[1](2)})(\a^{-2}(h)S^{-1}(\a^{-1}(a_{[1]})))_{[1]}]\b^{2}(a_{[0][1]})\\
&\stackrel{(4.1)}{=}& m_{[0]}\gamma(\a^{-2}(m_{[1]}))(\a^{-2}(h_{(1)})S^{-1}(\a^{-1}(a_{[1](2)})))\b(a_{[0][0]}) \\
&&\hspace{6cm}\o \a(\a^{-1}(h_{(2)})S^{-1}(a_{[1](1)}))\a^{2}(a_{[0][1]})\\
&=&m_{[0]}\gamma(\a^{-2}(m_{[1]}))(\a^{-2}(h_{(1)})S^{-1}(a_{[1](2)(2)}))a_{[0]}) \\
&&\hspace{6cm}\o (h_{(2)}S^{-1}(\a^{2}(a_{[1](2)(1)})))\a^{2}(a_{[1](1)})\\
&=&m_{[0]}\gamma(\a^{-2}(m_{[1]}))(\a^{-2}(h_{(1)})S^{-1}(\a^{-1}(a_{[1](2)})))a_{[0]}) \\
&&\hspace{6cm}\o \a(h_{(2)})(S^{-1}(\a^{2}(a_{[1](1)(2)}))\a^{2}(a_{[1](1)(1)}))\\
&=&m_{[0]}\gamma(\a^{-2}(m_{[1]}))(\a^{-2}(h_{(1)})S^{-1}(\a^{-2}(a_{[1]})))a_{[0]}) \o \a^{2}(h_{(2)}).
\end{eqnarray*}
Hence  $f$ is an epimorphism in $ \widetilde{\mathscr{H}}(\mathscr{M}_k)^{H}_{A}$ and has a $(H,\a)$-colinear section.
For any $k$-free presentation of $(M,\mu)$ in the category of $k$-modules
$
k^{I}\stackrel{\pi}{\rightarrow }M\rightarrow 0,
$
there is an epimorphism in $ \widetilde{\mathscr{H}}(\mathscr{M}_k)^{H}_{A}$:
\begin{eqnarray*}
(A|o H)^{I}\cong A\o H\o k \stackrel{g}{\rightarrow } M\rightarrow 0,
\end{eqnarray*}
where $g=f\circ (id_A\o id_H\o \pi)$.
Hence $A\o H$ is a generator in $ \widetilde{\mathscr{H}}(\mathscr{M}_k)^{H}_{A}$.
The proof is completed.$\hfill \Box$
\medskip

\section{The affineness criterion for relative Hom-Hopf modules}
\def\theequation{5. \arabic{equation}}
\setcounter{equation} {0} \hskip\parindent

{\bf Proposition 5.1.} Let $(H, \a)$ be a monoidal Hom-Hopf algebra with a bijective antipode and $(A,\b)$  a right
$(H,\a)$-Hom-comodule algebra.
 Assume that there exists a total quantum integral $\gamma: H\rightarrow Hom(H,A)$.
 Then $\widetilde{\rho}: A\rightarrow A\o H, \widetilde{\rho}(a)=a_{[0]}\o a_{[1]}$ splits in $\widetilde{\mathscr{H}}(\mathscr{M}_k)^{H}_{A}$.

{\bf Proof.} We consider the map
\begin{eqnarray*}
\lambda: A\o H\rightarrow A ~~~~~~~~~~\lambda(a\o h)=a_{[0]}\gamma(a_{[1]})(h)
\end{eqnarray*}
for all $a\in A$ and $h\in H$ is a right $(H,\a)$-colinear retraction of $\widetilde{\rho}$. In particular, by Eqs.(5.1) and 5.2, we obtain $\lambda(1_A\o 1_H)=1_A$, and
\begin{eqnarray}
\lambda(\beta^{-1}(a)\o h_{(1)})\o h_{(2)}=\lambda(a\o h)_{[0]}\o \lambda(a\o h)_{[1]}.
\end{eqnarray}
We define now
\begin{eqnarray*}
\Lambda: A\o H\rightarrow A ~~~~~~~~~~\Lambda(a\o h)=\lambda(1_A\o \a^{-1}(h)S^{-1}(a_{[1]}))\b(a_{[0]})
\end{eqnarray*}
for all $h\in H$ and $a\in A$. Then, for all $a\in A$ we have
\begin{eqnarray*}
(\Lambda\circ \widetilde{\rho} )(a)&=& \Lambda(a_{[0]}\o a_{[1]})\\
&=& \lambda(1_A\o \a^{-1}(a_{[1]})S^{-1}(a_{[0][1]}))\b(a_{[0][0]})\\
&=&\lambda(1_A\o a_{[1](2)}S^{-1}(a_{[1](1)}))a_{[0]}\\
&=& \lambda(1_A\o 1_H)\b^{-1}(a)=a,
\end{eqnarray*}
i.e., $\Lambda$ is still a retraction of $\widetilde{\rho}$. Now, for all $h\in H, a,b\in A$, we have
\begin{eqnarray*}
\Lambda((a\o h)b)&=&\Lambda(ab_{[0]}\o hb_{[1]})\\
&=&\lambda(1_A\o \a^{-1}(hb_{[1]})S^{-1}(a_{[1]}b_{[0][1]}))\b(a_{[0]}b_{[0][0]})\\
&=&\lambda(1_A\o h[\a^{-1}(b_{[1](2)})S^{-1}(\a^{-1}(b_{[1](1)}))]S^{-1}(a_{[1]}))\b(a_{[0]})b_{[0][0]}\\
&=&\lambda(1_A\o hS^{-1}(\a(a_{[1]})))[\b(a_{[0]})\b^{-1}(b)]\\
&=&\lambda(1_A\o \a^{-1}(h)S^{-1}(a_{[1]}))\b(a_{[0]})b\\
&=&\Lambda(a\o h)b
\end{eqnarray*}
Hence, $\Lambda$ is right $(A,\b)$-linear. It remains to prove that $\Lambda$ is also $(H, \a)$-colinear, we have
\begin{eqnarray*}
\widetilde{\rho}\Lambda(a\o h)&=&\widetilde{\rho}(\lambda(1_A\o \a^{-1}(h)S^{-1}(a_{[1]}))\b(a_{[0]}))\\
&=&\lambda(1_A\o \a^{-1}(h)S^{-1}(a_{[1]}))_{[0]}\b(a_{[0][0]})\o \lambda(1_A\o \a^{-1}(h)S^{-1}(a_{[1]}))_{[1]}\a(a_{[0][1]})\\
&\stackrel{(5.1)}{=}&\lambda(1_A\o \a^{-1}(h_{(1)})S^{-1}(a_{[1](2)}))\b(a_{[0][0]})\o [\a^{-1}(h_{(2))}S^{-1}(a_{[1](1)})]\a(a_{[0][1]})\\
&=&\lambda(1_A\o \a^{-1}(h_{(1)})S^{-1}(a_{[1](2)}))\b(a_{[0][0]})\o h_{(2)}[S^{-1}(a_{[1](1)})a_{[0][1]}]\\
&=&\lambda(1_A\o \a^{-1}(h_{(1)})S^{-1}(\a(a_{[1](2)(2)})))a_{[0]}\o h_{(2)}[S^{-1}(\a(a_{[1](2)(1)}))a_{[1](1)}]\\
&=&\lambda(1_A\o \a^{-1}(h_{(1)})S^{-1}(a_{[1](2)}))a_{[0]}\o h_{(2)}[S^{-1}(\a(a_{[1](1)(2)}))\a(a_{[1](1)(1)})]\\
&=&\lambda(1_A\o \a^{-1}(h_{(1)})S^{-1}(\a^{-1}(a_{[1]})))a_{[0]}\o \a(h_{(2)})\\
&=&(\Lambda\o id_H)\rho_{A\o H},
\end{eqnarray*}
i.e., proved that $\Lambda$ is a retraction of $\widetilde{\rho}$ in $\widetilde{\mathscr{H}}(\mathscr{M}_k)^{H}_{A}$. $\hfill \Box$

Recall from \cite{GC14} that the coinvariants of $A$ as
\begin{eqnarray*}
B=A^{coH}=\{a\in A|\widetilde{\rho}(a)=\b^{-1}(a)\o 1_H\}.
\end{eqnarray*}
Then $B$ is a subalgebra of $A$.

{\bf Proposition 5.2.} Let $(H, \a)$ be a monoidal Hom-Hopf algebra with a bijective antipode and $(A,\b)$  a right
$(H,\a)$-Hom-comodule algebra.  Assume that there exists $\gamma: H\rightarrow Hom(H,A)$ a total quantum integral. Then

(1) $B$ is a direct summand of $A$ as a left $B$-submodule,

(1) $B$ is a direct summand of $A$ as a right $B$-submodule.

{\bf Proof.} We shall prove that there exists a well defined left trace given by
the formula
\begin{eqnarray*}
t^{l}: (A, \b)\rightarrow (B, \b)~~~~~~~~~t^{l}(a)=\lambda(a\o 1_H)=a_{[0]}\gamma(a_{[1]})(1_H)
\end{eqnarray*}
for all $a\in A$, by Eq.(5.1) we obtain $\widetilde{\rho}(t^{l}(a))=t^{l}(a)\o 1_H$, ie., $t^{l}(a)\in B$, for all $a\in A$. Now for all $b\in B$
and $a\in A$, we have
\begin{eqnarray*}
t^{l}(ba)&=&b_{[0]}a_{[0]})\gamma(b_{[1]}a_{[1]})(1_H)\\
&=&[ba_{[0]}]\gamma(a_{[1]})(1_H)\\
&=& ba_{[0]}\gamma(a_{[1]})(1_H)\\
&=&bt^{l}(a),
\end{eqnarray*}
hence $t^{l}$ is a left $(B,\b)$-module map and finally
\begin{eqnarray*}
t^{l}(1_A)=1_A\gamma(1_H)(1_H)=1_A.
\end{eqnarray*}
hence $t^{l}$ is a left $(B,\b)$-module retraction of the inclusion $B\subset A$.

(2) Similarly, we can prove that the map given by the formula
\begin{eqnarray*}
t^{r}: (A, \b)\rightarrow (B, \b), ~~~~~~~~~t^{r}(a)=\Lambda(a\o 1_H)=\gamma(1_H)(S^{-1}(\a^2(a_{[1]})))\b(a_{[0]})
\end{eqnarray*}
for all $a\in A$, $t^{r}$ is a right $(B,\b)$-module retraction of the inclusion $B\subset A$. $\hfill \Box$

{\bf Definition 5.3.}  Let $(H, \a)$ be a monoidal Hom-Hopf algebra with a bijective antipode and $(A,\b)$  a right
$(H,\a)$-Hom-comodule algebra.  Assume that there exists $\gamma: H\rightarrow Hom(H,A)$ a total quantum integral. The map
\begin{eqnarray*}
t^l: (A, \b)\rightarrow (B, \b), ~~~~~~t^l(a)=a_{[0]}\gamma(a_{[1]})(1_H)
\end{eqnarray*}
for all $a\in A$ is called the quantum trace associated of $\gamma$.

Now, we shall construct functors connecting $\widetilde{\mathscr{H}}(\mathscr{M}_k)^{H}_{A}$ and $\widetilde{\mathscr{H}}(\mathscr{M}_k)_{B}$. First, if
$(M, \mu)\in \widetilde{\mathscr{H}}(\mathscr{M}_k)^{H}_{A}$, then
\begin{eqnarray*}
M^{coH}=\{m\in M|\widetilde{\rho}(m)=\mu^{-1}(m)\o 1_H\}.
\end{eqnarray*}
is the right $B$-module of the coinvariants of $M$. Furthermore, $M\rightarrow M^{coH}$ gives us a covariant functor
\begin{eqnarray*}
(-)^{coH}: \widetilde{\mathscr{H}}(\mathscr{M}_k)^{H}_{A}\rightarrow \widetilde{\mathscr{H}}(\mathscr{M}_k)_{B}.
\end{eqnarray*}

Now, for $(N,\nu)\in \widetilde{\mathscr{H}}(\mathscr{M}_k)_{B}, N\o_B A\in \widetilde{\mathscr{H}}(\mathscr{M}_k)^{H}_{A}$ via the structures
\begin{eqnarray*}
(a\o_B n)a'= a\b^{-1}(a')\o_B \nu(n)\\
\rho_{A\o_B N}(a\o_B n)=a_{[0]}\o \nu^{-1}(n)\o \a(a_{[1]})
\end{eqnarray*}
for all $n\in N, a,a'\in A$. In this way, we have constructed a covariant functor called the induction functor
\begin{eqnarray*}
A\o_B-: \widetilde{\mathscr{H}}(\mathscr{M}_k)_{B}\rightarrow \widetilde{\mathscr{H}}(\mathscr{M}_k)^{H}_{A}.
\end{eqnarray*}

{\bf Proposition 5.4.} Let $(H, \a)$ be a monoidal Hom-Hopf algebra with a bijective antipode and $(A,\b)$  a right
$(H,\a)$-Hom-comodule algebra. Then the induction functor $A\o_B-: \widetilde{\mathscr{H}}(\mathscr{M}_k)_{B}\rightarrow \widetilde{\mathscr{H}}(\mathscr{M}_k)^{H}_{A}$ is a left adjoint of the coinvariant functor $(-)^{coH}: \widetilde{\mathscr{H}}(\mathscr{M}_k)^{H}_{A}\rightarrow \widetilde{\mathscr{H}}(\mathscr{M}_k)_{B}$.

{\bf Proof.} Similar to \cite{MM02}.$\hfill \Box$

Recall from \cite{GC14} that $A\o H\in \widetilde{\mathscr{H}}(\mathscr{M}_k)^{H}_{A}$ with the structures
$$
\left\{
  \begin{array}{ll}
(a \otimes h ) \cdot b = a b _{[0]}\otimes h b _{[1]};\\
  \rho_{A\o H}  ( a \otimes h  )
  =(\b^{-1}(a) \otimes h _{(1)})\otimes \a(h _{(2)}),
\end{array}
\right.
$$
for any $a, b\in A$ and $h\in H$, and we identify $(A\o H)^{coH}\cong A $ via $a\o 1_H\mapsto a$. Then the the adjunction map
 can be viewed as a map in $\widetilde{\mathscr{H}}(\mathscr{M}_k)^{H}_{A}$, as
 \begin{eqnarray*}
 \psi: A\o_BA\rightarrow A\o H, ~~~~\psi(a\o b)=ab_{[0]}\o b_{[1]}
 \end{eqnarray*}
 for all $a, b\in A$. Here $A\o_B A\in \widetilde{\mathscr{H}}(\mathscr{M}_k)^{H}_{A}$ with the structures
 $$
\left\{
  \begin{array}{ll}
(a \otimes_B b ) \cdot a' = \b(a)\o_B b\b^{-1}(a');\\
  \rho_{A\o_B A}  ( a \otimes_B b  )
  =(\b^{-1}(a) \otimes b _{[0]})\otimes \a(b _{[1]}),
\end{array}
\right.
$$
 for all $a,a',b\in A$.

{\bf Definition 5.5.}  Let $(H, \a)$ be a monoidal Hom-Hopf algebra with a bijective antipode and $(A,\b)$  a right
$(H,\a)$-Hom-comodule algebra, and $B=A^{coH}$. Then $A/B$ is called a Galois extension if the the canonical map
\begin{eqnarray*}
\psi: A\o_BA\rightarrow A\o H, ~~~~\psi(a\o b)=\b^{-1}(a)b_{[0]}\o \a(b_{[1]})
\end{eqnarray*}
is bijective.

{\bf Theorem 5.6.} Let $(H, \a)$ be a monoidal Hom-Hopf algebra with a bijective antipode and $(A,\b)$  a right
$(H,\a)$-Hom-comodule algebra,  and $B=A^{coH} $. Assume that there exists $\gamma: H\rightarrow Hom(H,A)$ a total quantum integral.
Then
\begin{eqnarray*}
\eta_N: N\rightarrow (A\o_B N)^{coH},~~~~~~~\eta_N(n)=1_A\o_B n
\end{eqnarray*}
is an isomorphism of right $(B, \b)$-modules for all $(N,\nu)\in \widetilde{\mathscr{H}}(\mathscr{M}_k)_{B}$.

{\bf Proof.} Using the left quantum $t^l: (A, \b)\rightarrow (B, \b)$ we shall construct an inverse of $\eta_N$, we define
\begin{eqnarray*}
\theta_N: (A\o_B N)^{coH} \rightarrow N,~~~~~~~~~~\theta_N(a_i\o_B n_i)=\sum t^l(a_{i})n_i
\end{eqnarray*}
for any $a_i\o_B n_i\in (A\o_B N)^{coH}$. Since $t^l(1_A)=1_A$, and it is easy to check that $\theta_N\circ \eta_N=id_N$. Now, let $a_i\o_B n_i\in (A\o_B N)^{coH}$. Then
\begin{eqnarray*}
\a^{-1}(a_i)\o_B \nu^{-1}(n_i)\o 1_H=a_{[0]}\o_B \nu^{-1}(n_i)\o \a(a_{[1]}).
\end{eqnarray*}
It follows that
\begin{eqnarray*}
\a^{-1}(a_i)\o_B \nu^{-1}(n_i)\o 1_A=a_{[0]}\o_B \nu^{-1}(n_i)\o \g(a_{[1]})(1_H).
\end{eqnarray*}
Furthermore, we have
\begin{eqnarray*}
\a^{-1}(a_i) 1_A\o_B \nu^{-1}(n_i)=a_{[0]}\g(a_{[1]})(1_H)\o_B \nu^{-1}(n_i).
\end{eqnarray*}
then, we have
\begin{eqnarray*}
a_i\o_B n_i= t^l(a_{i})\o_B n_i.
\end{eqnarray*}
Hence, we have
\begin{eqnarray*}
(\eta_N\circ \theta_N)(a_i\o_B n_i)=\sum 1_A\o_Bt^l(a_{i})n_i=\sum t^l(a_{i})\o_B n_i=a_i\o_B n_i,
\end{eqnarray*}
Then $\theta_N$ is an inverse of $\eta_N$.$\hfill \Box$

We will now prove the main result of this section, that is, the affineness
criterion for relative Hom-Hopf modules.

{\bf Theorem 5.7. }   Let $(H, \a)$ be a monoidal Hom-Hopf algebra with a bijective antipode and $(A,\b)$  a right
$(H,\a)$-Hom-comodule algebra, and $B=A^{coH} $. Assume that

 (1) there exists a  total quantum integral $\gamma: H\rightarrow Hom(H,A)$;

 (2) the canonical map
$$
\beta: A\o_{B}A\rightarrow A\o H,\ \
a\o_{B}b\mapsto \b^{-1}(a)b_{[0]}\o \a(b_{[1]})
$$
is surjective. Then the induction functor $A\o_B-: \widetilde{\mathscr{H}}(\mathscr{M}_k)_{B}\rightarrow \widetilde{\mathscr{H}}(\mathscr{M}_k)^{H}_{A}$ is an equivalence of categories.

{\bf Proof. }
In Theorem 5.6 we have shown that the adjunction map $\eta_N: N\rightarrow (A\o_B N)^{coH}$
is an isomorphism for all $(N,\nu)\in \widetilde{\mathscr{H}}(\mathscr{M}_k)_{B}$. It remains to prove that the other
adjunction map, namely $\b_M: M^{coH}\o_B A\rightarrow M$ is also an isomorphism for all $(M,\mu) \in \widetilde{\mathscr{H}}(\mathscr{M}_k)^{H}_{A}$.

Let $(V, \omega)$ be a $k$-module. Then $A\o V\in \widetilde{\mathscr{H}}(\mathscr{M}_k)^{H}_{A}$ via the structures induced
 by $A$, i.e.,
\begin{eqnarray*}
(a\o v)a'= a\b^{-1}(a')\o \omega(v)\\
\rho_{A\o V}(a\o v)=a_{[0]}\o \omega^{-1}(v)\o \a(a_{[1]})
\end{eqnarray*}
for all $a, b\in A$ and $v\in V$.  In particular, for $V=A, A\o A\in \widetilde{\mathscr{H}}(\mathscr{M}_k)^{H}_{A}$ via
 \begin{eqnarray}
(a \otimes b ) \cdot a' =a\b^{-1}(a')\o \b(b);\\
  \rho_{A\o A}  ( a \otimes b  )
  =a_{[0]}\o \b^{-1}(b)\o \a(a_{[1]})
\end{eqnarray}
for all $a, b, a'\in A$. We will first prove that the adjunction map $\b_{A\o V}:(A\o V)^{coH}\o_BA\rightarrow A\o V$ is an isomorphism for any $k$-module $V$.

First, $V\o B$ and $B\o V \in \widetilde{\mathscr{H}}(\mathscr{M}_k)_{B}$ via the usual $B$-actions $(v\o a)\c b=\omega(v)\o a\b^{-1}(b)$, and $a'\c(b'\o v')=\b^{-1}(a')b'\o \omega(v')$ for all $a,b,a',b'\in B$ and $v,v'\in V$. The flip map $\tau: V\o B\rightarrow B\o V, ~~\tau(v\o b)=b\o v$, for all $b\in B$ and $v\in V$, is an isomorphism in  $\widetilde{\mathscr{H}}(\mathscr{M}_k)_{B}$. On the
other hand $A\o V \in \widetilde{\mathscr{H}}(\mathscr{M}_k)^{H}_{A}$ via the structures
induced by $A$, i.e.
 \begin{eqnarray*}
&&(v\o a)\c b= \omega(v)\o a\b^{-1}(b),\\
&&\r_{V\o A}(v\o a)= \omega^{-1}(v) \o a_{[0]} \o \a(a_{[1]}),
\end{eqnarray*}

It is easy to see that the flip map $\tau: A\o V\rightarrow V\o A, ~~\tau(a\o v)=v\o a$ is an isomorphism in $\widetilde{\mathscr{H}}(\mathscr{M}_k)^{H}_{A}$.

Applying Theorem 5.6 for $N=V\o B\cong B\o V$, we obtain the following isomorphisms in $\mathcal{M}_{B}$:
$$ B\o V\cong V\o B\cong(V\o B\o_B A)^{coH}\cong(V\o A)^{coH}\cong(A\o V)^{coH}.$$
 Hence, $(A\o V)^{coH}\o_BA\cong A\o V$.

Let us define
\begin{eqnarray*}
\widetilde{\psi}: A\o_{B}A\rightarrow A\o H, \ \ a\o_{B}b\mapsto \b^{-1}(a)b_{[0]}\o \a(b_{[1]}),
\end{eqnarray*}
for all $a,b \in A$. As $\psi$ is surjective, $\widetilde{\psi}$ is surjective, because
$\widetilde{\psi}=\psi\circ can$, where $can: A\o A\rightarrow A\o_B A$ is the canonical surjection.

Let us define now \begin{eqnarray*}
                  \xi: A\o A\rightarrow A\o H, ~~~\xi(a\o b)=(\widetilde{\psi}\circ \tau)(a\o b)=\b^{-1}(b)a_{[0]}\o \a(a_{[1]}),
                  \end{eqnarray*}
for all $a,b \in A$. The map $\xi$ is surjective, as $\widetilde{\psi}$ and $\tau$ are. We will prove that $\xi$ is a morphism in $\widetilde{\mathscr{H}}(\mathscr{M}_k)^{H}_{A}$. where $A\o A$ and $A\o H$ are realtive Hom-Hopf modules, respectively. Indeed,
\begin{eqnarray*}
\xi((a\o b)c)=\xi(a\b^{-1}(c)\o \b(b))&=&b[a_{(-1)}\b^{-1}(c_{(-1)})]\o \a(a_{(1)})(c_{(1)})\\
&=& [\b^{-1}(b)a_{(-1)}]c_{(-1)}\o \a(a_{(1)})(c_{(1)})\\
&=&\xi(a\o b)c,
\end{eqnarray*}
and
\begin{eqnarray*}
\rho_{A\o H}\xi(a \o b)&=&\rho_{A\o H}(\b^{-1}(b)a_{[0]}\o \a(a_{[1]}))\\
&=& [\b^{-2}(b)\b^{-1}(a_{[0]})\o \a(a_{[1](1)})]\o \a^{2}(a_{[1](2)})\\
&=& [\b^{-2}(b)a_{[0]}\o \a(a_{[0][1]})]\o \a(a_{[1]})\\
&=& (\xi \o id )(a_{[0]}\o \b^{-1}(b)\o \b^{-1}(a_{[1]}))\\
&=& ( \xi\o id_H)\rho_{A\o A}(a\o b).
\end{eqnarray*}
Hence, $\xi$ is a surjective morphism in $\widetilde{\mathscr{H}}(\mathscr{M}_k)^{H}_{A}$. $A\o H$ is projective as a left $A$-module,
where $A\o H$ is a left $A$-module in the usual way. On the other hand, by Proposition 3.1 we obtain the map
\begin{eqnarray*}
u: A\o H\rightarrow A\o H, \ \ a\o h\mapsto \b(a_{[0]})\o h\a^{-1}(a_{[1]})
\end{eqnarray*}
is an isomorphism of left $A$-modules.            It follows that there exits $\zeta: A\o H \rightarrow A\o A$ such that
$\xi\circ\zeta=id_{H\boxtimes A}$ since $A\o A\rightarrow A\o H$ is surjective. Hence, $\xi$ splits
in the category of right $A$-modules. In particular $\xi$ is a $k$-split epimorphism in $\widetilde{\mathscr{H}}(\mathscr{M}_k)^{H}_{A}$.

Let now $(M,\mu)\in \widetilde{\mathscr{H}}(\mathscr{M}_k)^{H}_{A}$. Then $A\o A\o M\in \widetilde{\mathscr{H}}(\mathscr{M}_k)^{H}_{A}$ via the structures arising from $A\o A$, that is,
\begin{eqnarray}
&&(a\o b\o m)\c c=a\b^{-2}(c)\o \b(b) \o \mu(m);\\
&&\rho_{A\o A\o M}(a\o b\o m)=\b^{-1}(a)\o b_{[0]}\o \mu^{-1}(m) \o \a^{2}(b_{[1]}),
\end{eqnarray}
for all $a, b,c \in A$ and $m\in M$. On the other hand, $A\o H\o M\in \widetilde{\mathscr{H}}(\mathscr{M}_k)^{H}_{A}$ via the structures arising from $A\o H$, that is,
 \begin{eqnarray*}
(a\o h\o m)b= a\b^{-1}(b_{[0]})\o h\a^{-1}(b_{[1]})\o \mu(m)\\
\rho_{A\o H\o M}(a\o h\o m)= \b^{-1}(a)\o h_{(1)}\o \mu^{-1}(m) \o \a^{2}(h_{(2)}).
 \end{eqnarray*}
for all $a, b \in A, h\in H$ and $m\in M$. We obtain that
\begin{eqnarray*}
\xi\o id_M: A\o A\o M\rightarrow A\o H\o M
\end{eqnarray*}
is a $k$-split epimorphism in $\widetilde{\mathscr{H}}(\mathscr{M}_k)^{H}_{A}$.

 By Theorem 4.8 we obtain that the
map
\begin{eqnarray*}
f: A\o H\o M\rightarrow M, \ \
a\o h\o m\mapsto \mu(m_{[0]})\gamma(\a^{-1}(m_{[1]}))(\a^{-2}(h)S^{-1}(\a^{-1}(a_{[1]})))\b(a_{[0]})
\end{eqnarray*}
is a $k$-split epimorphism in $\widetilde{\mathscr{H}}(\mathscr{M}_k)^{H}_{A}$. Hence, the composition
$$
g=f\circ(\xi\o id_M): A\o A\o M\rightarrow M,
$$
$$
a\o b\o m\mapsto \mu(m_{[0]})\gamma(\a^{-1}(m_{[1]}))(S^{-1}(b_{[1]}))b_{[0]}\b^{-1}(a)
$$
is a $k$-split epimorphism in $\widetilde{\mathscr{H}}(\mathscr{M}_k)^{H}_{A}$. We note that the structure of $A\o A
\o M$ as an object in $\widetilde{\mathscr{H}}(\mathscr{M}_k)^{H}_{A}$ is of the form $A\o V$, for the $k$-module
$V=A\o M$.

To conclude, we have constructed a $k$-split epimorphism in $\widetilde{\mathscr{H}}(\mathscr{M}_k)^{H}_{A}$
\begin{eqnarray*}
A\o A\o M=(M_1,\pi)\stackrel{g}{\longrightarrow} (M,\mu)\longrightarrow 0
\end{eqnarray*}
such that the adjunction map $\psi_{M_1}$ for $(M_1,\pi)$ is bijective. As $g$ is $k$-split and there exists a  total quantum integral $\gamma: H\rightarrow Hom(H,A)$, we obtain that
$g$ also splits in $\widetilde{\mathscr{H}}(\mathscr{M}_k)^{H}$. In particular, the sequence
\begin{eqnarray*}
(M_1^{coH}, \pi)\stackrel{g^{coH}}{\longrightarrow} (M^{coH},\mu)\longrightarrow 0
\end{eqnarray*}
 is exact. Continuing the resolution with $Ker (g)$ instead of $M$, we obtain an
exact sequence in $\widetilde{\mathscr{H}}(\mathscr{M}_k)^{H}_{A}$
\begin{eqnarray*}
(M_2,P)\longrightarrow (M_1,\pi)\longrightarrow (M,\mu)\longrightarrow0
\end{eqnarray*}
which splits in $\widetilde{\mathscr{H}}(\mathscr{M}_k)^{H}$ and the adjunction maps for $(M_1,\pi)$ and $(M_2,P)$ are bijective.
Using the Five lemma we obtain that the adjunction map for $(M,\mu)$ is
bijective.$\hfill \Box$

Finally we consider a special case. In this case if setting $A=H$, then $(A,\b)$ is a right $H$-comodule algebra in
a natural way, and we can define the coinvariants of $(H,\a)$ as
\begin{eqnarray*}
B=H^{coH}=\{h\in H|\tilde{\r}(h)=\a^{-1}(h)\o 1_H \}
\end{eqnarray*}
then $(B,\a)$ is a subalgebra of $(H,\a)$. Hence we can obtain the following result.

{\bf Corollary 5.8. }  Let $(H, \a)$ be a monoidal Hom-Hopf algebra with a bijective antipode and $(A,\b)$  a right
$(H,\a)$-Hom-comodule algebra, and $B=H^{coH}$. Assume that
 \\
 (1) there exists a  total quantum integral $\gamma: H\rightarrow Hom(H,H)$;\\
 (2) the canonical map
$$
\psi: H\o_{B}H\rightarrow H\o H,\ \
h\o_{B}g\mapsto \a^{-1}(h)g_{(1)}\o \a(g_{(2)})
$$
is surjective. Then the induction functor $-\o_{B}H:\widetilde{\mathscr{H}}(\mathscr{M}_k)_{B}\rightarrow \widetilde{\mathscr{H}}(\mathscr{M}_k)^{H}_{H}$ is an equivalence of categories.

\section*{Acknowledgements}

 The work is supported by the NSF of Jiangsu Province (No. BK2012736) and  the Fund of Science and Technology Department of Guizhou Province (No. 2014GZ81365).

\end{document}